\documentclass{article}

\usepackage{arxiv}
\usepackage{mypackage}

\usepackage[utf8]{inputenc} 
\usepackage[T1]{fontenc}    
\usepackage{hyperref}       
\usepackage{url}            
\usepackage{booktabs}       
\usepackage{amsfonts}       
\usepackage{nicefrac}       
\usepackage{microtype}      
\usepackage{lipsum}

\makeatletter
\def\and{%
  \end{tabular}%
  \begin{tabular}[t]{c}}%
\def\@fnsymbol#1{\ensuremath{\ifcase#1\or a\or b\or c\or
   d\or e\or f\or g\or h\or i\else\@ctrerr\fi}}
\makeatother

\begin{document}
\title{A Dynamic Subspace Based BFGS Method for Large Scale Optimization Problem}

\author{
    Zheng Li,
    Shi Shu,
    Jian-Ping Zhang
}
\maketitle

\begin{abstract}
Large-scale unconstrained optimization is a fundamental and important class of, yet not well-solved problems in numerical optimization. The main challenge in designing an algorithm is to require a few storage locations or very inexpensive computations while preserving global convergence. In this work, we propose a novel approach solving large-scale unconstrained optimization problem by combining the dynamic subspace technique and the BFGS update algorithm.
It is clearly demonstrated that our approach has the same rate of convergence in the dynamic subspace as the BFGS and less memory than L-BFGS. Further, we give the convergence analysis by constructing the mapping of low-dimensional Euclidean space to the adaptive subspace. We compare our hybrid algorithm with the BFGS and L-BFGS approaches. Experimental results show that our hybrid algorithm offers several significant advantages such as parallel computing, convergence efficiency, and robustness.%
\blfootnote{
    Zheng Li, Shi Shu, Jian-Ping Zhang \\
    School of Mathematics and Computational Science,Xiangtan University,Xiangtan, Hunan, 411105, China \\
    E-mail: lizheng.math.ai@gmail.com \\
    E-mail: shushi@xtu.edu.cn \\
    E-mail: jpzhang@xtu.edu.cn
}
\end{abstract}

\section{Introduction}
\label{sec: Introduction}
In the past two decades, the applications of subspace optimization in various specific problems have been extensively studied. The main aim of designing a conjugate gradient-like method on a subspace for large-scale unconstrained optimization problems is to reduce the overall cost of computation and computational storage.
For example, Narkiss et al. in \cite{narkiss2005sequential} proposed a sequential subspace optimization (SESOP) method. At each iteration, the search for a minimum of the objective function over a subspace spanned by the current gradient direction and a few previous steps. For convex problems, the method orders the rate of convergence to be $\frac{1}{N^2}$ ($N$ is the number of the iterations). Andrei in \cite{andrei2014accelerated} given a three-term conjugate gradient algorithm for large-scale unconstrained optimization using subspace minimizing technique. Similar to the SECOP method, the subspace of this method is also spanned by several specific vectors. The numerical experiments show that this new algorithm is more robust than conjugate gradient algorithms respectively proposed by Hestenes and Stiefel \cite{hestenes1952methods}, Dai and Liao \cite{dai2001new}, Dai and Yuan and Polak \cite{dai1999nonlinear}, Ribiére and Poliak \cite{polak1969note}, as well as the limited memory quasi-Newton method (L-BFGS method reported in \cite{gilbert1989some}) and the discrete truncated-Newton method (TN method described in \cite{nash1991numerical} by Nash ).

These subspace-based methods have a significant improvement for solving large-scale optimization problems, but they are rarely based on the quasi-Newton theory. 
However,  Wang and Yuan in \cite{wang2006subspace} proposed some subspace trust-region algorithms by studying subspace properties of trust-region methods for unconstrained optimization, In fact, every iteration of Wang and Yuan method requires only $O(n^2)$ floating
point operations and a gradient evaluation. limited-memory are not required during the process.
To fix this gap, we construct a highly efficient subspace (Fast-BFGS) method for large-scale smooth unconstrained problems.

Moreover, we observe from the BFGS method proposed by Fletcher et al.\cite{fletcher1987practical} in solving large-scale optimization problems,  if the inverse Hessian approximation $H_k$ can be expressed as the truncated form (\ref{eqn: truncated form}) with a special initial matrix,
the computational cost associated with the search direction $-H_k \nabla f_k$ and the storage will be greatly reduced. 
Unfortunately, $H_k$ cannot be equivalently expressed as the truncated form (\ref{eqn: truncated form}) for general large-scale unconstrained optimization problems. One of our main ideas is to consider how to equivalently truncate $H_k$. 

Our motivation employing the subspace method is that, if $H_k$ is replaced by a specially constructed matrix $\tilde{H}_k$ in our giving truncated form (\ref{eqn: truncated form}), thus the proposed iteration $\boldsymbol{x}_k$ of the minimum $\boldsymbol{x} \in \mathbb{R}^n$ will be constrained on a low-dimensional hyperplane $\mathcal{P}$ (the dimension of $\mathcal{P}$ is denoted by $m$). Moreover, by constructing a linear mapping from $\mathbb{R}^m$ to $\mathcal{P}$, we can show that the updating process of $\boldsymbol{x} \in \mathcal{P}$ is equivalent to the variable $\boldsymbol{\xi} \in \mathbb{R}^m$ update in $\boldsymbol{\xi}$-minimization applying the standard BFGS method. The dynamic hyperplane $\mathcal{P}_k$ will gradually stabilize and closer to the minimum of $f$ when $\| \nabla f_k \|_2 \rightarrow 0$. Thus, the Fast-BFGS method has the same rate of convergence as the BFGS method near the minimum point. Our proposed Fast-BFGS method has the following advantages:

\begin{itemize}
\item We propose
a novel dynamic subspace-based BFGS method solving large-scale unconstrained minimization problem, where the search direction at each iteration is constrained on a low-dimensional dynamic subspace and updated by using the BFGS method. 
\item Our method only needs to store $m$ n-dimensional vectors and a $m \times m$ matrix. Many experiments show that it's still effective when $m$ tends to be small, therefore we have an advantage over the L-BFGS method (described by Liu et al.\cite{liu1989limited}) in limited memory, and also inexpensive cost in parallel computing.
\item The numerical results on a variety of problems from the CUTE collection\cite{bongartz1995cute,lukvsan2010modified,andrei2008unconstrained} show that our algorithm is more effective than the BFGS method or L-BFGS method in the vast majority of cases.
\end{itemize}

The present paper is organized as follows. 
Section \ref{sec: BFGS Update Formula in Truncated Form} proves that $H_k$ in the traditional BFGS method can be decomposed into the sum of rank-one matrices with a group of special direction vector $\{\tilde{\boldsymbol{s}}_i\}_{i < k}$, then gives the necessary and sufficient conditions for $H_k$ to be expressed as the truncated form (\ref{eqn: truncated form}).
A new Hessian-free method is presented in Section \ref{sec: Fast-BFGS for Large-Scale Optimization}. 
Convergence analysis of our method is given in Section \ref{sec: Convergence Analysis}, and it is shown that our method is actually equivalent to the standard BFGS method on a low-dimensional subspace.
Section \ref{sec: Modify the Search Direction} combines a method of adaptively updating the subspace with the method described in Section \ref{sec: Fast-BFGS for Large-Scale Optimization} to enable our algorithm to solve global optimization problems.
Finally, Section \ref{sec: Experiments} presents numerical results on a variety of problem in the CUTE collection.

\section{The proposing truncated BFGS update}
\label{sec: BFGS Update Formula in Truncated Form}
In this section, we provide some notations about BFGS method and our dynamic subspace which are useful for our later analysis.
For any $k \ge 0$, the variations of variables and gradients are denoted by 
\begin{equation*}
\boldsymbol{s}_k = \boldsymbol{x}_{k+1} - \boldsymbol{x}_k, \quad \boldsymbol{y}_k = \nabla f_{k+1} - \nabla f_k.
\end{equation*}
The rescaling of $\boldsymbol{s}_k$ and $\boldsymbol{y}_k$ are denoted by
\begin{equation*}
\tilde{\boldsymbol{s}}_k = \frac{\boldsymbol{s}_k}{\sqrt{\boldsymbol{s}_k^T \boldsymbol{y}_k}}, \quad \tilde{\boldsymbol{y}}_k = \frac{\boldsymbol{y}_k}{\sqrt{\boldsymbol{s}_k^T \boldsymbol{y}_k}}.
\end{equation*}
The block matrix of stacks $\begin{bmatrix} \tilde{\boldsymbol{s}}_0 & \cdots & \tilde{\boldsymbol{s}}_{k-1} \end{bmatrix}$ is denoted by $S_k$. 
For any given $m < n$, a matrix $H \in \mathbb{R}^{n \times n}$ is called in a truncated form if there exists $k > 0$ and $L \in \mathbb{R}^{\min(k, m) \times \min(k, m)}$, make $H$ can be decomposed as 
\begin{equation}
H = \tilde{S}_k L \tilde{S}_k^T.
\label{eqn: truncated form}
\end{equation} 
where
\begin{equation*}
\tilde{S}_k =  
\begin{cases}
\begin{bmatrix} \tilde{\boldsymbol{s}}_0 & \cdots & \tilde{\boldsymbol{s}}_{k-1} \end{bmatrix}, & \forall k < m; \\
\begin{bmatrix} \tilde{\boldsymbol{s}}_{k-m} & \cdots & \tilde{\boldsymbol{s}}_{k-1} \end{bmatrix}, & \forall k \ge m.
\end{cases}
\end{equation*}
When $k>m$, the function $T_k: \mathbb{R}^m \rightarrow \mathbb{R}^{m \times m}$ is given by
\begin{equation}
T_k(\boldsymbol{t}) = \begin{bmatrix} t_1 & 1 & \cdots & 0 \\ \vdots & \vdots & \ddots & \vdots \\ t_{m-1} & 0 & \cdots & 1 \\ t_m - \tilde{\boldsymbol{y}}_k^T\tilde{\boldsymbol{s}}_{k-m} & -\tilde{\boldsymbol{y}}_k^T\tilde{\boldsymbol{s}}_{k+1-m} & \cdots & -\tilde{\boldsymbol{y}}_k^T \tilde{\boldsymbol{s}}_{k-1} \end{bmatrix},
\label{eqn: BFGS tk function}
\end{equation}
where $\boldsymbol{t} = \begin{bmatrix} t_1 & \cdots & t_m \end{bmatrix}^T$.




Now assume that the BFGS method is start with $m$ step and 
\begin{equation*}
H_m = \sum_{i=0}^{m-1} \sum_{j=0}^{m-1} (L_{m-1})_{ij} \tilde{\boldsymbol{s}}_i \tilde{\boldsymbol{s}}_j^T = S_m L_{m-1} S_m^T.
\end{equation*}
We obtain
\begin{align*}
    (\mathrm{BFGS}) \quad H_{m+1} &= (I - \tilde{\boldsymbol{s}}_m \tilde{\boldsymbol{y}}_m^T) H_m (I - \tilde{\boldsymbol{s}}_m \tilde{\boldsymbol{y}}_m^T)^T + \tilde{\boldsymbol{s}}_m \tilde{\boldsymbol{s}}_m^T \notag \\
    &= (I - \tilde{\boldsymbol{s}}_m \tilde{\boldsymbol{y}}_m^T) 
    S_m L_{m-1} S_m^T
    (I - \tilde{\boldsymbol{s}}_m \tilde{\boldsymbol{y}}_m^T)^T + \tilde{\boldsymbol{s}}_m \tilde{\boldsymbol{s}}_m^T \notag \\
    &= (I - \tilde{\boldsymbol{s}}_m \tilde{\boldsymbol{y}}_m^T) 
    \begin{bmatrix}
    \tilde{\boldsymbol{s}}_0 & \cdots & \tilde{\boldsymbol{s}}_{m-1}
    \end{bmatrix}
    L_{m-1}
    \begin{bmatrix}
    \tilde{\boldsymbol{s}}_0^T \\ \vdots \\ \tilde{\boldsymbol{s}}_{m-1}^T
    \end{bmatrix} 
    (I - \tilde{\boldsymbol{s}}_m \tilde{\boldsymbol{y}}_m^T)^T + \tilde{\boldsymbol{s}}_m \tilde{\boldsymbol{s}}_m^T \notag \\
    &= \begin{bmatrix}
    \tilde{\boldsymbol{s}}_0 & \dots & \tilde{\boldsymbol{s}}_m
    \end{bmatrix}
    \begin{bmatrix} I \\ - \tilde{\boldsymbol{y}}_m^T S_m \end{bmatrix} L_{m-1} \begin{bmatrix} I & -S_m^T \tilde{\boldsymbol{y}}_m \end{bmatrix}
    \begin{bmatrix}
    \tilde{\boldsymbol{s}}_0^T \\ \vdots \\ \tilde{\boldsymbol{s}}_m^T
    \end{bmatrix} + \tilde{\boldsymbol{s}}_m \tilde{\boldsymbol{s}}_m^T \notag \\
    &= \begin{bmatrix}
    \tilde{\boldsymbol{s}}_0 & \dots & \tilde{\boldsymbol{s}}_m
    \end{bmatrix}
    \left( \begin{bmatrix} I \\ - \tilde{\boldsymbol{y}}_m^T S_m \end{bmatrix} L_{m-1} \begin{bmatrix} I & -S_m^T \tilde{\boldsymbol{y}}_m \end{bmatrix} + 
    \begin{bmatrix} \boldsymbol{0} & \boldsymbol{0} \\ \boldsymbol{0}^T & 1 \\ \end{bmatrix} 
    \right)
    \begin{bmatrix}
    \tilde{\boldsymbol{s}}_0^T \\ \vdots \\ \tilde{\boldsymbol{s}}_m^T
    \end{bmatrix}.
\end{align*}

A easy induction gives
\begin{equation}
H_k = S_k L_{k-1} S_k^T = \sum_{i=1}^{k} \sum_{j=1}^{k} (L_{k-1})_{ij} \tilde{\boldsymbol{s}}_i \tilde{\boldsymbol{s}}_j^T, \quad \forall k \ge m,
\label{eqn: BFGS hk}
\end{equation}
where
\begin{equation}
L_k = \begin{bmatrix} I \\ - \tilde{\boldsymbol{y}}_k^T S_k \end{bmatrix} L_{k-1} \begin{bmatrix} I & -S_k^T \tilde{\boldsymbol{y}}_k \end{bmatrix} + \begin{bmatrix} \boldsymbol{0} & \boldsymbol{0} \\ \boldsymbol{0}^T & 1 \\ \end{bmatrix}, \quad \forall k \ge m.
\label{eqn: ak}
\end{equation}

\begin{theorem}
\label{thm: spd ak}
If $L_{m-1}$ is symmetrical positive definite, the matrix $L_k$ is symmetrical positive definite for $k \ge m$.
\end{theorem}
\begin{proof}
The Eq. (\ref{eqn: ak}) follows that if $L_{k-1}$ is symmetrical semi-positive definite then $L_k$ is also symmetrical semi-positive definite. Decomposing the matrix $L_k$ into 
\begin{equation*}
L_k = \begin{bmatrix} I \\ - \tilde{\boldsymbol{y}}_k^T S_k \end{bmatrix} L_{k-1} \begin{bmatrix} I & -S_k^T \tilde{\boldsymbol{y}}_k \end{bmatrix} + \begin{bmatrix} \boldsymbol{0} & \boldsymbol{0} \\ \boldsymbol{0}^T & 1 \end{bmatrix} = \begin{bmatrix} L_{k-1}  & \boldsymbol{0} \\ -\tilde{\boldsymbol{y}}_k^T S_k L_{k-1} & 1 \end{bmatrix} \begin{bmatrix} I  & S_k^T \tilde{\boldsymbol{y}}_k \\ \boldsymbol{0}^T & 1 \end{bmatrix}^{-1}
\end{equation*}
we have if $\mathrm{det}(L_{k-1}) \neq 0$ then $\mathrm{det}(L_k) \neq 0$ is also not equal
to $0$. Continuing by induction, the proof is completed.
\end{proof}

This theorem reveals the following important phenomenon.

\begin{theorem}
\label{thm: rank hk}
Since
\begin{equation*}
H_k \nabla f_k = S_k L_k S_k^T \nabla f_k = \sum_{i=1}^k (L_k S_k^T \nabla f_k)_i \tilde{\boldsymbol{s}}_{i-1} \in \mathrm{span} \{\tilde{\boldsymbol{s}}_0, \cdots, \tilde{\boldsymbol{s}}_{k-1} \}, \quad \forall k \ge m,
\end{equation*}
we have 
\begin{equation*}
\tilde{\boldsymbol{s}}_k \in \mathrm{span} \{\tilde{\boldsymbol{s}}_0, \cdots, \tilde{\boldsymbol{s}}_{k-1} \}, \quad \forall k \ge m.
\end{equation*}
It follows that
\begin{equation*}
\tilde{\boldsymbol{s}}_k \in \mathrm{span} \{\tilde{\boldsymbol{s}}_0, \cdots, \tilde{\boldsymbol{s}}_{m-1} \}, \mathrm{rank} (S_k) = \mathrm{rank} (S_m), \quad \forall k \ge m.
\end{equation*} 
If $\mathrm{rank} (L_{m-1}) = m$, then the proof of Theorem (\ref{thm: spd ak}) shows that $\mathrm{rank} (L_{k-1}) = k$. Thus, we finally have 
\begin{equation}
\mathrm{rank} (H_k) = m \Leftrightarrow \mathrm{rank} (H_m) = m.
\end{equation}
\end{theorem}

Consider the large-scale unconstrained optimization problems. The Eq. (\ref{eqn: BFGS hk}) shows that the search direction in BFGS method can be computed by using 
\begin{equation}
\boldsymbol{p}_k = -H_k \nabla f_k = -\begin{bmatrix} \tilde{\boldsymbol{s}}_0 & \cdots & \tilde{\boldsymbol{s}}_{k-1} \end{bmatrix} \left( L_{k-1} \left( \begin{bmatrix} \tilde{\boldsymbol{s}}_0^T \\ \vdots \\ \tilde{\boldsymbol{s}}_{k-1}^T \end{bmatrix} \nabla f_k \right) \right), \quad \forall k \ge k_0.
\label{eqn: pk}
\end{equation}
Since the computational complexity of $\boldsymbol{p}_k$ is grows as $k$ increases, that is not good for designing Eq. (\ref{eqn: pk}) as a parallel program. However, if there exist a integer $m \ll n$ and a series of $m \times m$ matrices $\{ \tilde{L}_{k-1} \}_{k \ge m}$ satisfy 
\begin{equation}
H_k = \begin{bmatrix} \tilde{\boldsymbol{s}}_{k-m} & \cdots & \tilde{\boldsymbol{s}}_{k-1} \end{bmatrix} \tilde{L}_{k-1} \begin{bmatrix} \tilde{\boldsymbol{s}}_{k-m}^T \\ \cdots \\ \tilde{\boldsymbol{s}}_{k-1}^T \end{bmatrix}, \quad \forall k \ge m,
\label{eqn: BFGS hk limited form}
\end{equation}
then the computational complexity of $\boldsymbol{p}_k$ is always equal to $O(mn) + O(m^2)$. 
\begin{remark}
The update formula of L-BFGS method is also in a truncated form. Different from Eq. (\ref{eqn: BFGS hk limited form}), The L-BFGS method makes an approximate estimation of $H_k \nabla f_k$, and needs the data of $\{ \boldsymbol{y}_{k-m}, \cdots, \boldsymbol{y}_{k-1} \}$ when calculating $\boldsymbol{p}_k$.
\end{remark}

The following theorem will shows us that the $H_k$ in the BFGS method cannot be expressed  in the truncated form Eq. (\ref{eqn: BFGS hk limited form}) for general optimization problems.

\begin{theorem}
\label{thm: BFGS hk truncated form}
Let $m$ be an integer subject to $m \le n$. If $\mathrm{rank}(\tilde{S}_m) = m$, then for $\forall k \ge m$ and any objective function, there exist a $m \times m$ matrix $\hat{L}_{k-1}$ subject to
\begin{equation*}
H_k = \tilde{S}_k \hat{L}_{k-1} \tilde{S}_k^T,
\end{equation*}
if and only if 
\begin{equation}
\mathrm{rank}(\tilde{S}_k) \equiv m, \quad k \ge m. 
\label{eqn: BFGS hk limited form condition}
\end{equation}
\end{theorem}

\begin{proof}
We first show the sufficiency.

When $k=m$, we have $H_m = S_m L_{m-1} S_m^T = \tilde{S}_m L_{m-1} \tilde{S}_m^T$. 

When $k>m$, assuming there exist a series matrices $\{ \hat{L}_m, \cdots, \hat{L}_{k-1} \}$ subject to 
\begin{equation*}
H_i = \tilde{S}_i \hat{L}_{i-1} \tilde{S}_i^T, \quad i=m, \cdots, k-1,
\end{equation*}
then we have 
\begin{equation*}
H_i \nabla f_i = \tilde{S}_i (\hat{L}_{i-1} \tilde{S}_i^T \nabla f_i) = \sum_{j=1}^{m} (\hat{L}_{i-1} \tilde{S}_i^T \nabla f_i)_j \tilde{\boldsymbol{s}}_{i-1-m+j} \in \mathrm{span} \{\tilde{\boldsymbol{s}}_{i-m}, \cdots, \tilde{\boldsymbol{s}}_{i-1} \}, \quad i=m, \cdots, k-1.
\end{equation*}
Using this equation and noting $H_i \nabla f_i$ is parallel to $\tilde{\boldsymbol{s}}_i$,it lead to 
\begin{equation*}
\tilde{\boldsymbol{s}}_0, \cdots, \tilde{\boldsymbol{s}}_k \in \mathrm{span\ } \{\tilde{s}_0, \cdots, \tilde{s}_{m-1} \}.
\end{equation*}
Applying the conclusion to the condition Eq. (\ref{eqn: BFGS hk limited form condition}), there exist a $m$ dimension vector $\boldsymbol{t}^{(k)}$ which satisfy
\begin{equation}
\tilde{\boldsymbol{s}}_{k-m} = \tilde{S}_{k+1} \boldsymbol{t}^{(k)}.
\label{eqn: subquestion-0}
\end{equation} 
Notice
\begin{align*}
(I - \tilde{\boldsymbol{s}}_k \tilde{\boldsymbol{y}}_k^T) \tilde{S}_k 
&= \tilde{S}_{k+1} T_k(\boldsymbol{0}) + \begin{bmatrix} \tilde{\boldsymbol{s}}_{k-m} & \boldsymbol{0}_{n \times m-1} \end{bmatrix} \\
&= \tilde{S}_{k+1} T_k(\boldsymbol{0}) + \tilde{S}_{k+1} \begin{bmatrix} \boldsymbol{t}^{(k)} & \boldsymbol{0}_{n \times m-1} \end{bmatrix} \\
&= \tilde{S}_{k+1} T_k(\boldsymbol{t}^{(k)}), 
\end{align*}
we have
\begin{align*}
(\mathrm{BFGS}) \quad H_{k+1} &= (I - \tilde{\boldsymbol{s}}_k \tilde{\boldsymbol{y}}_k^T) H_k (I - \tilde{\boldsymbol{y}}_k \tilde{\boldsymbol{s}}_k^T) + \tilde{\boldsymbol{s}}_k \tilde{\boldsymbol{s}}_k^T \\
&= (I - \tilde{\boldsymbol{s}}_k \tilde{\boldsymbol{y}}_k^T) \tilde{S}_k \hat{L}_{k-1} \tilde{S}_k^T (I - \tilde{\boldsymbol{y}}_k \tilde{\boldsymbol{s}}_k^T) + \tilde{\boldsymbol{s}}_k \tilde{\boldsymbol{s}}_k^T \\
&= \tilde{S}_{k+1} T_k(\boldsymbol{t}^{(k)}) \hat{L}_{k-1} (T_k(\boldsymbol{t}^{(k)}))^T \tilde{S}_{k+1}^T + \tilde{\boldsymbol{s}}_k \tilde{\boldsymbol{s}}_k^T \\
&= \tilde{S}_{k+1} \hat{L}_k \tilde{S}_{k+1}^T, 
\end{align*}
where 
\begin{equation*}
\hat{L}_k = T_k(\boldsymbol{t}^{(k)}) \hat{L}_{k-1} (T_k(\boldsymbol{t}^{(k)}))^T + \begin{bmatrix} \boldsymbol{0} & \boldsymbol{0} \\ \boldsymbol{0}^T & 1 \\ \end{bmatrix}, 
\end{equation*}
Continuing by induction, the proof of sufficiency is completed.

We are now truning to the proof of necessity.

For $\forall k \ge m$, notice
\begin{equation*}
H_k \nabla f_k = \tilde{S}_k \tilde{L}_{k-1} \tilde{S}_k^T \nabla f_k = \sum_{i=1}^{m} (\tilde{L}_{k-1} \tilde{S}_k^T \nabla f_k)_i \tilde{\boldsymbol{s}}_{k+i-m-1} \in \mathrm{span} \{ \tilde{\boldsymbol{s}}_{k-m}, \cdots, \tilde{\boldsymbol{s}}_{k-1} \},
\end{equation*}
thus $\tilde{\boldsymbol{s}}_{k} \in \mathrm{span} \{ \tilde{\boldsymbol{s}}_{k-m}, \cdots, \tilde{\boldsymbol{s}}_{k-1} \}$. It follows that
\begin{equation}
\mathrm{rank} (\tilde{S}_k) \le \mathrm{rank} (\tilde{S}_{k-1}), \quad k \ge m.
\label{eqn: BFGS hk limited form proof-2}
\end{equation}

For any objective function $f$, noting
\begin{align*}
(\mathrm{BFGS}) \quad H_{k+1} \nabla f_{k+1} &= (I - \tilde{\boldsymbol{s}}_k \tilde{\boldsymbol{y}}_k^T) H_k (I - \tilde{\boldsymbol{y}}_k \tilde{\boldsymbol{s}}_k^T) \nabla f_{k+1} + (\tilde{\boldsymbol{s}}_k^T \nabla f_{k+1}) \tilde{\boldsymbol{s}}_k \\
&= (I - \tilde{\boldsymbol{s}}_k \tilde{\boldsymbol{y}}_k^T) \tilde{S}_k \left( \hat{L}_{k-1} \tilde{S}_k^T (I - \tilde{\boldsymbol{y}}_k \tilde{\boldsymbol{s}}_k^T) \nabla f_{k+1} \right) + (\tilde{\boldsymbol{s}}_k^T \nabla f_{k+1}) \tilde{\boldsymbol{s}}_k \\
&= \tilde{S}_{k+1} \left( T_k(\boldsymbol{0}) \hat{L}_{k-1} \tilde{S}_k^T (I - \tilde{\boldsymbol{y}}_k \tilde{\boldsymbol{s}}_k^T) \nabla f_{k+1} \right) \\
&+ \begin{bmatrix} \tilde{\boldsymbol{s}}_{k-m} & \boldsymbol{0}_{n \times m-1} \end{bmatrix} \left( \hat{L}_{k-1} \tilde{S}_k^T (I - \tilde{\boldsymbol{y}}_k \tilde{\boldsymbol{s}}_k^T) \nabla f_{k+1} \right) + (\tilde{\boldsymbol{s}}_k^T \nabla f_{k+1}) \tilde{\boldsymbol{s}}_k \\
&= \sum_{i=1}^{m} \left( T_k(\boldsymbol{0}) \hat{L}_{k-1} \hat{L}_{k-1} \tilde{S}_k^T (I - \tilde{\boldsymbol{y}}_k \tilde{\boldsymbol{s}}_k^T) \nabla f_{k+1} \right)_i \tilde{\boldsymbol{s}}_{k+i-m} + (\tilde{\boldsymbol{s}}_k^T \nabla f_{k+1}) \tilde{\boldsymbol{s}}_k \\
&+ \left( \hat{L}_{k-1} \tilde{S}_k^T (I - \tilde{\boldsymbol{y}}_k \tilde{\boldsymbol{s}}_k^T) \nabla f_{k+1} \right)_1 \tilde{\boldsymbol{s}}_{k-m} \\
&\in \left( \hat{L}_{k-1} \tilde{S}_k^T (I - \tilde{\boldsymbol{y}}_k \tilde{\boldsymbol{s}}_k^T) \nabla f_{k+1} \right)_1 \tilde{\boldsymbol{s}}_{k-m} + \mathrm{span} \{ \tilde{\boldsymbol{s}}_{k+1-m}, \cdots, \tilde{\boldsymbol{s}}_{k} \}.
\end{align*}
and
\begin{equation*}
H_{k+1} \nabla f_{k+1} = \tilde{S}_{k+1} \hat{L}_k \tilde{S}_{k+1}^T \nabla f_{k+1} = \sum_{i=1}^{m} (\hat{L}_k \tilde{S}_{k+1}^T \nabla f_{k+1})_i \tilde{\boldsymbol{s}}_{k+i-m} \in \mathrm{span} \{ \tilde{\boldsymbol{s}}_{k+1-m}, \cdots, \tilde{\boldsymbol{s}}_k \},
\end{equation*}
we have 
\begin{equation}
\tilde{\boldsymbol{s}}_{k-m} \in \mathrm{span} \{ \tilde{\boldsymbol{s}}_{k+1-m}, \cdots, \tilde{\boldsymbol{s}}_{k} \}.
\label{eqn: BFGS hk limited form proof-3}
\end{equation}
Using Eq. (\ref{eqn: BFGS hk limited form proof-3}) in to Eq. (\ref{eqn: BFGS hk limited form proof-2}), we finally have 
\begin{equation*}
\mathrm{rank} (\tilde{S}_k) = \mathrm{rank} (\tilde{S}_{k-1}).
\end{equation*}
Continuing by induction, the proof of necessity is completed.
\end{proof}

From this theorem, 
The condition (\ref{eqn: BFGS hk limited form condition}) is too difficult to satisfy so that we cannot use the truncated form (\ref{eqn: BFGS hk limited form}) to calculate the search direction most of the time.
 In the next section, we will modify the calculation of $\boldsymbol{t}^{(k)}$ and and propose a new truncated form.

\section{Fast-BFGS for Large-Scale Optimization}
\label{sec: Fast-BFGS for Large-Scale Optimization}

\subsection {Fast-BFGS Updating}
We proposed the Fast-BFGS method which is based on the following formula 
\begin{equation*}
(\mathrm{Fast-BFGS}) \quad \tilde{H}_k = \tilde{S}_k \tilde{L}_{k-1} \tilde{S}_k^T, \quad \tilde{L}_k = T_k(\boldsymbol{t}^{(k)}) \tilde{L}_{k-1} (T_k(\boldsymbol{t}^{(k)}))^T + \begin{bmatrix} \boldsymbol{0} & \boldsymbol{0} \\ \boldsymbol{0}^T & 1 \\ \end{bmatrix}, \quad \forall k \ge m,
\end{equation*}
where $\boldsymbol{t}^{(k)}$ is computed by the following unconstrained problem
\begin{equation}
\boldsymbol{t}^{(k)} = \mathop{\arg\min}_{\boldsymbol{t} \in \mathbb{R}^m} \| \tilde{S}_{k+1} \boldsymbol{t} - \tilde{\boldsymbol{s}}_{k-m} \|_2.
\label{eqn: subquestion-1}
\end{equation}
An argument similar to the proof of Theorem \ref{thm: BFGS hk truncated form} shows that the variable $\boldsymbol{x}$ iterated by 
\begin{equation}
\boldsymbol{x}_{k+1} =\boldsymbol{x}_k - \tau_k \tilde{H}_k  \nabla f_k, \tau_k = \mathop{\arg\min}_{\tau} f(\boldsymbol{x}_k - \tau \tilde{H}_k \nabla f_k), \quad k \ge m,
\label{eqn: fast BFGS xk}
\end{equation}
converges to the solution of the problem
\begin{equation*}
\min f(\boldsymbol{x}), \quad \mathrm{subject\ to\ } \quad \boldsymbol{x} \in \boldsymbol{x}_0 + \bigcap_{k=m}^{\infty} \mathrm{span} \{ \tilde{\boldsymbol{s}}_{k-m}, \cdots, \tilde{\boldsymbol{s}}_{k-1} \}.
\end{equation*}

\subsection {Relationship with BFGS Method}
\label{sec: Relationship with BFGS Method}

Let $(\boldsymbol{\beta}_0, \cdots, \boldsymbol{\beta}_m)$ be the reordering of $(\tilde{\boldsymbol{s}}_{k-m}, \cdots, \tilde{\boldsymbol{s}}_k)$ and satisfy 
\begin{equation*}
\boldsymbol{\beta}_0 \in \mathrm{span} \{ \boldsymbol{\beta}_1, \cdots, \boldsymbol{\beta}_m \}, \boldsymbol{\beta}_m = \tilde{\boldsymbol{s}_k},
\end{equation*}
If we want $\{ \boldsymbol{x}_k \}_{k \ge m}$ converges to the solution of the constrained optimization problem
\begin{equation}
\min f(\boldsymbol{x}), \quad \mathrm{subject\ to\ } \quad \boldsymbol{x} \in \boldsymbol{x}_0 + \mathrm{span} \{ \tilde{\boldsymbol{s}}_0, \cdots, \tilde{\boldsymbol{s}}_{m-1} \}
\label{eqn: constrained optimization problem}
\end{equation}
with the rate of convergence is superlinear, we need to set $\mathrm{rank} (\tilde{H}_m) = m$, then, for $k \ge m$, reset
\begin{equation*}
\tilde{\boldsymbol{s}}_{k-m} \leftarrow \boldsymbol{\beta}_0, \cdots, \tilde{\boldsymbol{s}}_{k-1} \leftarrow \boldsymbol{\beta}_{m-1}
\end{equation*}
after updating $\boldsymbol{x}_k$ to $\boldsymbol{x}_{k+1}$.

We will prove the proposition in the next section.

\section{Convergence Analysis} 
\label{sec: Convergence Analysis}

In this section, for $\forall k \ge m$, we always reorder $(\tilde{\boldsymbol{s}}_{k-m}, \cdots, \tilde{\boldsymbol{s}}_k)$ after updating $\boldsymbol{x}_k$ to $\boldsymbol{x}_{k+1}$ which is described in subsection \ref{sec: Relationship with BFGS Method}. 

Assuming $\mathrm{rank} (S_m) = m$ and denoting the Schmidt orthogonalization result of $S_m = \begin{bmatrix} \tilde{\boldsymbol{s}}_0, \cdots, \tilde{\boldsymbol{s}}_{m-1} \end{bmatrix}$ as $S_m^{\mathrm{unit}}$, then consider to solve the equivalent form of problem \ref{eqn: constrained optimization problem}
\begin{equation*}
\mathop{\min}_{\boldsymbol{\xi} \in \mathbb{R}^m} f(\boldsymbol{x}_0 + S_m^{\mathrm{unit}} \boldsymbol{\xi})
\end{equation*}
by using BFGS method 
\begin{align}
\boldsymbol{\xi}_{k+1} &= \boldsymbol{\xi}_k - \tau_k^{\boldsymbol{\xi}} H_k^{\boldsymbol{\xi}} \nabla_{\boldsymbol{\xi}} f_k \label{eqn: fast BFGS superlinear-1} \\
(\mathrm{BFGS}) \quad H_{k+1}^{\boldsymbol{\xi}} &= \left( I_m - \tilde{\boldsymbol{s}}_k^{\boldsymbol{\xi}} (\tilde{\boldsymbol{y}}_k^{\boldsymbol{\xi}})^T \right) H_k^{\boldsymbol{\xi}} \left( I_m - \tilde{\boldsymbol{y}}_k^{\boldsymbol{\xi}} (\tilde{\boldsymbol{s}}_k^{\boldsymbol{\xi}})^T \right) + \tilde{\boldsymbol{s}}_k^{\boldsymbol{\xi}} (\tilde{\boldsymbol{s}}_k^{\boldsymbol{\xi}})^T \label{eqn: fast BFGS superlinear-2}
\end{align}
where 
\begin{align*}
\tau_k^{\boldsymbol{\xi}} &= \mathop{\arg\min}_{\tau} f(\boldsymbol{x}_0 + S_m^{\mathrm{unit}} (\boldsymbol{\xi}_k - \tau H_k^{\boldsymbol{\xi}} \nabla_{\boldsymbol{\xi}} f_k)) \\
\nabla_{\boldsymbol{\xi}} f_k &= \frac{\partial f(\boldsymbol{x}_0 + S_m^{\mathrm{unit}} \boldsymbol{\xi})}{\partial \boldsymbol{\xi}} \bigg|_{\boldsymbol{\xi} = \boldsymbol{\xi}_k}, \\
\boldsymbol{s}_k^{\boldsymbol{\xi}} &= \boldsymbol{\xi}_{k+1} - \boldsymbol{\xi}_{k}, \\
\boldsymbol{y}_k^{\boldsymbol{\xi}} &= \nabla_{\boldsymbol{\xi}} f_{k+1} - \nabla_{\boldsymbol{\xi}} f_k, \\
\tilde{\boldsymbol{s}}_k^{\boldsymbol{\xi}} &= \frac{\boldsymbol{s}_k^{\boldsymbol{\xi}}}{\sqrt{(\boldsymbol{s}_k^{\boldsymbol{\xi}})^T \boldsymbol{y}_k^{\boldsymbol{\xi}}}}, \\
\tilde{\boldsymbol{y}}_k^{\boldsymbol{\xi}} &= \frac{\boldsymbol{y}_k^{\boldsymbol{\xi}}}{\sqrt{(\boldsymbol{s}_k^{\boldsymbol{\xi}})^T \boldsymbol{y}_k^{\boldsymbol{\xi}}}}.
\end{align*}
The iteration of Eq. (\ref{eqn: fast BFGS superlinear-1}) and Eq. (\ref{eqn: fast BFGS superlinear-2}) starts with $\boldsymbol{\xi} = \boldsymbol{\xi}_{m-1}$ and $H_{m}^{\boldsymbol{\xi}} = (S_m^{\mathrm{unit}})^T \tilde{H}_m S_m^{\mathrm{unit}}$ is a matrix of full rank.

The following lemma shows the relationship of $\{ \boldsymbol{x}_k \}_{k \ge m}$ and $\{ \boldsymbol{\xi}_k \}_{k \ge m}$.

\begin{lemma}
\label{lem: hk transform}
For $\forall k \ge m$, let $\boldsymbol{x}_k = \boldsymbol{x}_0 + S_m^{\mathrm{unit}} \boldsymbol{\xi}_k$, then the equation 
\begin{equation}
\begin{cases}
\tau_k^{\boldsymbol{\xi}} &= \mathop{\arg\min}_{\tau} f(\boldsymbol{x}_k - \tau \tilde{H}_k \nabla f_k); \\
H_{k}^{\boldsymbol{\xi}} &= (S_m^{\mathrm{unit}})^T \tilde{H}_k S_m^{\mathrm{unit}}
\end{cases}
\label{eqn: hk transform}
\end{equation}
holds.
\end{lemma}

\begin{proof}
We divide our proof in four steps.

First, we need to verify two properties. For $\forall k \ge m$, notice $\boldsymbol{s}_k = S_m^{\mathrm{unit}} \boldsymbol{s}_k^{\boldsymbol{\xi}}, \boldsymbol{y}_k^{\boldsymbol{\xi}} = (S_m^{\mathrm{unit}})^T \boldsymbol{y}_k$ and $(S_m^{\mathrm{unit}})^T S_m^{\mathrm{unit}} = I_m$, we have
\begin{equation}
I_m - \tilde{\boldsymbol{s}}_k^{\boldsymbol{\xi}} (\tilde{\boldsymbol{y}}_k^{\boldsymbol{\xi}})^T = (S_m^{\mathrm{unit}})^T (I_n - \tilde{\boldsymbol{s}}_k \tilde{\boldsymbol{y}}_k^T) S_m^{\mathrm{unit}}, \quad k \ge m.
\label{eqn: hk transform-1}
\end{equation}
For $\forall \boldsymbol{s} \in \mathrm{span} \{ \tilde{\boldsymbol{s}}_0, \cdots, \tilde{\boldsymbol{s}}_{m-1} \}$, we have the inversion formula
\begin{equation*}
S_m^{\mathrm{unit}} (S_m^{\mathrm{unit}})^T \boldsymbol{s} = \boldsymbol{s}.
\end{equation*}
It is follows that
\begin{equation}
S_m^{\mathrm{unit}} (S_m^{\mathrm{unit}})^T \tilde{S}_k = \tilde{S}_k, \quad k \ge 1.
\label{eqn: hk transform-2}
\end{equation}

The next thing to do in the proof is to show 
\begin{equation*}
\tau_m^{\boldsymbol{\xi}} = \mathop{\arg\min}_{\tau} f(\boldsymbol{x}_m - \tau \tilde{H}_m \nabla f_m).
\end{equation*}
Notice $H_{m}^{\boldsymbol{\xi}} = (S_m^{\mathrm{unit}})^T \tilde{H}_m S_m^{\mathrm{unit}}$, we obtain 
\begin{align*}
\tau_m^{\boldsymbol{\xi}} 
&= \mathop{\arg\min}_{\tau} f(\boldsymbol{x}_0 + S_m^{\mathrm{unit}} (\boldsymbol{\xi}_m - \tau H_m^{\boldsymbol{\xi}} \nabla_{\boldsymbol{\xi}} f_m)) \\
&= \mathop{\arg\min}_{\tau} f(\boldsymbol{x}_m - \tau S_m^{\mathrm{unit}} H_m^{\boldsymbol{\xi}} (S_m^{\mathrm{unit}})^T \nabla f_m) \\
&= \mathop{\arg\min}_{\tau} f(\boldsymbol{x}_m - \tau \tilde{H}_m \nabla f_m).
\end{align*}

Another step in the proof is to assume 
there exist $k > m$ makes $(\tau_{m+1}^{\boldsymbol{\xi}}, H_{m+1}^{\boldsymbol{\xi}}), \cdots, (\tau_k^{\boldsymbol{\xi}}, H_k^{\boldsymbol{\xi}})$ satisfy Eq. (\ref{eqn: hk transform}). 

Finally, we have to show that the Eq. (\ref{eqn: hk transform}) still holds when $k \leftarrow k+1$. Combining Eq. (\ref{eqn: hk transform-1}) and Eq. (\ref{eqn: hk transform-2}), and then using $\tilde{H}_k = \tilde{S}_k \tilde{L}_{k-1} \tilde{S}_k^T$, we get
\begin{align}
H_{k+1}^{\boldsymbol{\xi}}
=& \left( I_m - \tilde{\boldsymbol{s}}_k^{\boldsymbol{\xi}} (\tilde{\boldsymbol{y}}_k^{\boldsymbol{\xi}})^T \right) H_k^{\boldsymbol{\xi}} \left( I_m - \tilde{\boldsymbol{y}}_k^{\boldsymbol{\xi}} (\tilde{\boldsymbol{s}}_k^{\boldsymbol{\xi}})^T \right) + \tilde{\boldsymbol{s}}_k^{\boldsymbol{\xi}} (\tilde{\boldsymbol{s}}_k^{\boldsymbol{\xi}})^T \notag \\
=& (S_m^{\mathrm{unit}})^T \left( (I_n - \tilde{\boldsymbol{s}}_k \tilde{\boldsymbol{y}}_k^T) S_m^{\mathrm{unit}} H_k^{\boldsymbol{\xi}} (S_m^{\mathrm{unit}})^T (I_n - \tilde{\boldsymbol{y}}_k \tilde{\boldsymbol{s}}_k^T)\right) S_m^{\mathrm{unit}} + \tilde{\boldsymbol{s}}_k^{\boldsymbol{\xi}} (\tilde{\boldsymbol{s}}_k^{\boldsymbol{\xi}})^T \notag \\
=& (S_m^{\mathrm{unit}})^T \left( (I_n - \tilde{\boldsymbol{s}}_k \tilde{\boldsymbol{y}}_k^T) \left( S_m^{\mathrm{unit}} (S_m^{\mathrm{unit}})^T \tilde{S}_k \right) \tilde{L}_{k-1} \left( S_m^{\mathrm{unit}} (S_m^{\mathrm{unit}})^T \tilde{S}_k \right)^T (I_n - \tilde{\boldsymbol{y}}_k \tilde{\boldsymbol{s}}_k^T)\right) S_m^{\mathrm{unit}} + \tilde{\boldsymbol{s}}_k^{\boldsymbol{\xi}} (\tilde{\boldsymbol{s}}_k^{\boldsymbol{\xi}})^T \notag \\
=& (S_m^{\mathrm{unit}})^T \left( (I_n - \tilde{\boldsymbol{s}}_k \tilde{\boldsymbol{y}}_k^T) \tilde{S}_k \tilde{L}_{k-1} \tilde{S}_k^T (I_n - \tilde{\boldsymbol{y}}_k \tilde{\boldsymbol{s}}_k^T)\right) S_m^{\mathrm{unit}} + \tilde{\boldsymbol{s}}_k^{\boldsymbol{\xi}} (\tilde{\boldsymbol{s}}_k^{\boldsymbol{\xi}})^T.
\label{eqn: hk transform-3}
\end{align}
According to the definition of the reordering process which is described in subsection \ref{sec: Relationship with BFGS Method}, $\tilde{\boldsymbol{s}}_{k-m} \in \mathrm{span} \{ \tilde{\boldsymbol{s}}_{k-m+1}, \cdots, \tilde{\boldsymbol{s}}_k \}$, it follows that $\mathop{\min}_{\boldsymbol{t} \in \mathbb{R}^m} \| \tilde{S}_{k+1} \boldsymbol{t} - \tilde{\boldsymbol{s}}_{k-m} \|_2 = 0$. Using the notation $T_k(\cdot)$ with is defined by Eq. (\ref{eqn: BFGS tk function}), we have 
\begin{align}
& (S_m^{\mathrm{unit}})^T (I - \tilde{\boldsymbol{s}}_k \tilde{\boldsymbol{y}}_k^T) \tilde{S}_k \notag \\
=& (S_m^{\mathrm{unit}})^T \tilde{S}_{k+1} T_k(\boldsymbol{0}) + \begin{bmatrix} (S_m^{\mathrm{unit}})^T \tilde{\boldsymbol{s}}_{k-m} & \boldsymbol{0}_{n \times m-1} \end{bmatrix} \notag \\
=& (S_m^{\mathrm{unit}})^T \tilde{S}_{k+1}  T_k(\boldsymbol{t}^{(k)}) + \begin{bmatrix} (S_m^{\mathrm{unit}})^T (\tilde{\boldsymbol{s}}_{k-m} - \tilde{S}_{k+1} \boldsymbol{t}^{(k)}) & \boldsymbol{0}_{n \times m-1} \end{bmatrix} \notag \\
=& (S_m^{\mathrm{unit}})^T \tilde{S}_{k+1} T_k(\boldsymbol{t}^{(k)}).
\label{eqn: hk transform-4}
\end{align}
Applying Eq. (\ref{eqn: hk transform-4}) to Eq. (\ref{eqn: hk transform-3}), we have, in light of $\boldsymbol{s}_k = S_m^{\mathrm{unit}} \boldsymbol{s}_k^{\boldsymbol{\xi}}, \boldsymbol{y}_k^{\boldsymbol{\xi}} = (S_m^{\mathrm{unit}})^T \boldsymbol{y}_k$ and $(S_m^{\mathrm{unit}})^T S_m^{\mathrm{unit}} = I_m$ 
\begin{align*}
H_{k+1}^{\boldsymbol{\xi}}
=& (S_m^{\mathrm{unit}})^T \left( \tilde{S}_{k+1} T_k(\boldsymbol{t}^{(k)}) \tilde{L}_{k-1} (T_k(\boldsymbol{t}^{(k)}))^T \tilde{S}_{k+1}^T \right) S_m^{\mathrm{unit}} + (S_m^{\mathrm{unit}})^T \tilde{\boldsymbol{s}}_k \tilde{\boldsymbol{s}}_k^T S_m^{\mathrm{unit}} \notag \\
=& (S_m^{\mathrm{unit}})^T \left( \tilde{S}_{k+1} \left( T_k(\boldsymbol{t}^{(k)}) \tilde{L}_{k-1} (T_k(\boldsymbol{t}^{(k)}))^T + \begin{bmatrix} \boldsymbol{0} & \boldsymbol{0} \\ \boldsymbol{0}^T & 1 \end{bmatrix} \right) \tilde{S}_{k+1}^T \right) S_m^{\mathrm{unit}} \notag \\
=& (S_m^{\mathrm{unit}})^T \left( \tilde{S}_{k+1} \tilde{L}_k \tilde{S}_{k+1}^T \right) S_m^{\mathrm{unit}} \notag \\
=& (S_m^{\mathrm{unit}})^T \tilde{H}_{k+1} S_m^{\mathrm{unit}}.
\end{align*}
then
\begin{align*}
\tau_{k+1}^{\boldsymbol{\xi}} 
&= \mathop{\arg\min}_{\tau} f(\boldsymbol{x}_0 + S_m^{\mathrm{unit}} (\boldsymbol{\xi}_{k+1} - \tau H_{k+1}^{\boldsymbol{\xi}} \nabla_{\boldsymbol{\xi}} f_{k+1})) \\
&= \mathop{\arg\min}_{\tau} f(\boldsymbol{x}_{k+1} - \tau S_m^{\mathrm{unit}} H_{k+1}^{\boldsymbol{\xi}} (S_m^{\mathrm{unit}})^T \nabla f_{k+1}) \\
&= \mathop{\arg\min}_{\tau} f(\boldsymbol{x}_{k+1} - \tau \tilde{H}_{k+1} \nabla f_{k+1}).
\end{align*}
Continuing by induction, the proof is computed.
\end{proof}

\begin{corollary}
From Theorem \ref{thm: rank hk} and Lemma \ref{lem: hk transform}, for $k \ge m$, $\mathrm{rank} \tilde{H}_k = m$ holds with $\mathrm{rank} \tilde{H}_m = m$.
\end{corollary}

Using Lemma \ref{lem: hk transform}, we can prove that $\{ \boldsymbol{x}_k \}_{k \ge m}$ is convergent in the rate of convergence is superlinear.

\begin{theorem}
\label{thm: superlinear}
If the $\boldsymbol{t}^{(k)}$ is defined by problem \ref{eqn: subquestion-1} and $\mathrm{rank} (S_m) = m$, then the sequence $\{ \boldsymbol{x}_k \}_{k \ge m}$ generated by Eq. (\ref{eqn: fast BFGS xk}) converges to the solution of constrained problem \ref{eqn: constrained optimization problem} with the rate of convergence is superlinear.
\end{theorem}

\begin{proof}
From Lemma \ref{lem: hk transform} and Eq. (\ref{eqn: hk transform-2}), we have
\begin{align}
\boldsymbol{x}_{k+1} - \boldsymbol{x}_k &= -\tau_k \tilde{H}_k \nabla f_k = -\tau_k^{\boldsymbol{\xi}} \tilde{H}_k \nabla f_k = -\tau_k^{\boldsymbol{\xi}} \tilde{S}_k \tilde{L}_{k-1} \tilde{S}_k^T \nabla f_k \notag \\
&= -\tau_k^{\boldsymbol{\xi}} \left( S_m^{\mathrm{unit}} (S_m^{\mathrm{unit}})^T \tilde{S}_k \right) \tilde{L}_{k-1} \left( S_m^{\mathrm{unit}} (S_m^{\mathrm{unit}})^T \tilde{S}_k \right)^T \nabla f_k \notag \\
&= -\tau_k^{\boldsymbol{\xi}} S_m^{\mathrm{unit}} \left( (S_m^{\mathrm{unit}})^T \tilde{H}_k S_m^{\mathrm{unit}} \right) \left( (S_m^{\mathrm{unit}})^T \nabla f_k \right) \notag \\
&= -\tau_k^{\boldsymbol{\xi}} S_m^{\mathrm{unit}} H_k^{\boldsymbol{\xi}} \nabla_{\boldsymbol{\xi}} f_k \notag \\
&= S_m^{\mathrm{unit}} (\boldsymbol{\xi}_{k+1} - \boldsymbol{\xi}_k).
\label{eqn: superlinear-1}
\end{align}
Notice $(S_m^{\mathrm{unit}})^T S_m^{\mathrm{unit}} = I_m$, we have
\begin{equation}
\|\boldsymbol{x}_k - \boldsymbol{x}_{k-1}\|_2 = \|\boldsymbol{\xi}_k - \boldsymbol{\xi}_{k-1}\|_2 
\end{equation}
in light of $\{ \boldsymbol{\xi}_k \}_{k \ge m}$ is generated by BFGS method, the proof will be completed.
\end{proof}

Lemma \ref{lem: hk transform} also shows that $\tilde{H}_k$ satisfy the secant equation.

\begin{theorem}
\label{thm: secant equation}
For $k \ge m$,
\begin{equation*}
\tilde{H}_{k+1} \boldsymbol{y}_k = \boldsymbol{s}_k.
\end{equation*}
\end{theorem}

\begin{proof}
Noting the $H_k^{\boldsymbol{\xi}}$ is updated by the BFGS method, it follows that
\begin{equation}
\tilde{\boldsymbol{s}}_k^{\boldsymbol{\xi}}
= H_{k+1}^{\boldsymbol{\xi}} \tilde{\boldsymbol{y}}_k^{\boldsymbol{\xi}}
, \quad k \ge m.
\label{eqn: secant equation proof-1}
\end{equation}
Applying Lemma \ref{lem: hk transform} to Eq. (\ref{eqn: secant equation proof-1}) and noting $\boldsymbol{y}_k^{\boldsymbol{\xi}} = (S_m^{\mathrm{unit}})^T \boldsymbol{y}_k$, we have 
\begin{equation}
\tilde{\boldsymbol{s}}_k^{\boldsymbol{\xi}}
= (S_m^{\mathrm{unit}})^T \tilde{H}_{k+1} S_m^{\mathrm{unit}} (S_m^{\mathrm{unit}})^T \tilde{\boldsymbol{y}}_k 
, \quad k \ge m.
\label{eqn: secant equation proof-2}
\end{equation}
Expanding $\tilde{H}_{k+1}$ into $\tilde{S}_{k+1} \tilde{L}_k \tilde{S}_{k+1}^T$, then applying the inversion formula $S_m^{\mathrm{unit}} (S_m^{\mathrm{unit}})^T \tilde{S}_{k+1} = \tilde{S}_{k+1}$ to Eq. (\ref{eqn: secant equation proof-2}), we obtain
\begin{equation}
S_m^{\mathrm{unit}} \tilde{\boldsymbol{s}}_k^{\boldsymbol{\xi}}
= \left( S_m^{\mathrm{unit}} (S_m^{\mathrm{unit}})^T \tilde{S}_{k+1} \right) \tilde{L}_k \left( S_m^{\mathrm{unit}} (S_m^{\mathrm{unit}})^T \tilde{S}_{k+1} \right)^T \tilde{\boldsymbol{y}}_k
=\tilde{S}_{k+1} \tilde{L}_k \tilde{S}_{k+1}^T \tilde{\boldsymbol{y}}_k
= \tilde{H}_{k+1} \tilde{\boldsymbol{y}}_k
, \quad k \ge m.
\label{eqn: secant equation proof-3}
\end{equation}
From the inversion formula $\tilde{\boldsymbol{s}}_k = S_m^{\mathrm{unit}} \left( (S_m^{\mathrm{unit}})^T \tilde{\boldsymbol{s}}_k \right)$ and $\tilde{\boldsymbol{s}}_k = S_m^{\mathrm{unit}} \tilde{\boldsymbol{s}}_k^{\boldsymbol{\xi}}$, we have
\begin{equation}
\tilde{\boldsymbol{s}}_k 
= S_m^{\mathrm{unit}} \left( (S_m^{\mathrm{unit}})^T \tilde{\boldsymbol{s}}_k \right) 
= S_m^{\mathrm{unit}} \tilde{\boldsymbol{s}}_k^{\boldsymbol{\xi}}
, \quad k \ge m.
\label{eqn: secant equation proof-4}
\end{equation}
Combining (\ref{eqn: secant equation proof-4}) and (\ref{eqn: secant equation proof-3}), finally, we have
\begin{align*}
\boldsymbol{s}_k = \sqrt{\boldsymbol{s}_k^T \boldsymbol{y}_k} \tilde{\boldsymbol{s}}_k
= \sqrt{\boldsymbol{s}_k^T \boldsymbol{y}_k} S_m^{\mathrm{unit}} \tilde{\boldsymbol{s}}_k^{\boldsymbol{\xi}}  
= \sqrt{\boldsymbol{s}_k^T \boldsymbol{y}_k} \tilde{H}_{k+1} \tilde{\boldsymbol{y}}_k   
= \tilde{H}_{k+1} \boldsymbol{y}_k, \quad k \ge m.
\end{align*}
The prove is completed.
\end{proof}

\section{Modify the Search Direction} 
\label{sec: Modify the Search Direction}

In last section, Theorem \ref{thm: superlinear} shows that Fast-BFGS is used in solving the constrained problem \ref{eqn: constrained optimization problem}. In order to arrive at a global minimum point, we modify the search direction as 
\begin{equation*}
\boldsymbol{p}_k = -\tilde{H}_k \nabla f_k - \boldsymbol{v}_k, \quad k \ge 0
\end{equation*}
where the definition of $\tilde{H}_k = \tilde{S}_k \tilde{L}_{k-1} \tilde{S}_k$ is extended as 
\begin{equation}
\begin{cases}
\tilde{L}_0 &= \begin{bmatrix} 1 \end{bmatrix} \\
\tilde{L}_k &= \begin{bmatrix} I \\ - \tilde{\boldsymbol{y}}_k^T S_k \end{bmatrix} \tilde{L}_{k-1} \begin{bmatrix} I & -  S_k^T \tilde{\boldsymbol{y}}_k \end{bmatrix} + \begin{bmatrix} \boldsymbol{0} & \boldsymbol{0} \\ \boldsymbol{0}^T & 1 \\ \end{bmatrix}, \quad \forall 1 \le k < m, \\
\tilde{L}_k &= T_k(\boldsymbol{t}^{(k)}) \tilde{L}_{k-1} (T_k(\boldsymbol{t}^{(k)}))^T + \begin{bmatrix} \boldsymbol{0} & \boldsymbol{0} \\ \boldsymbol{0}^T & 1 \\ \end{bmatrix}, \quad \forall k \ge m,
\end{cases}
\label{eqn: Fast-BFGS}
\end{equation}

Consider to make the rate of convergence to be superlinear, we need to find the suitable $\boldsymbol{v}_k$ by
\begin{equation}
\boldsymbol{v}_k = \mathop{\arg\min}_{\boldsymbol{v}^T \nabla f_k > 0} \| \nabla f(\boldsymbol{x}_k - \tilde{H}_k \nabla f_k - \boldsymbol{v}) \|_2.
\label{eqn: unsolve problem}
\end{equation}
The crucial basis in our modification of search direction is the following figure.
\begin{figure}[H]
\centering
\includegraphics[width=1.0\textwidth, angle=0]{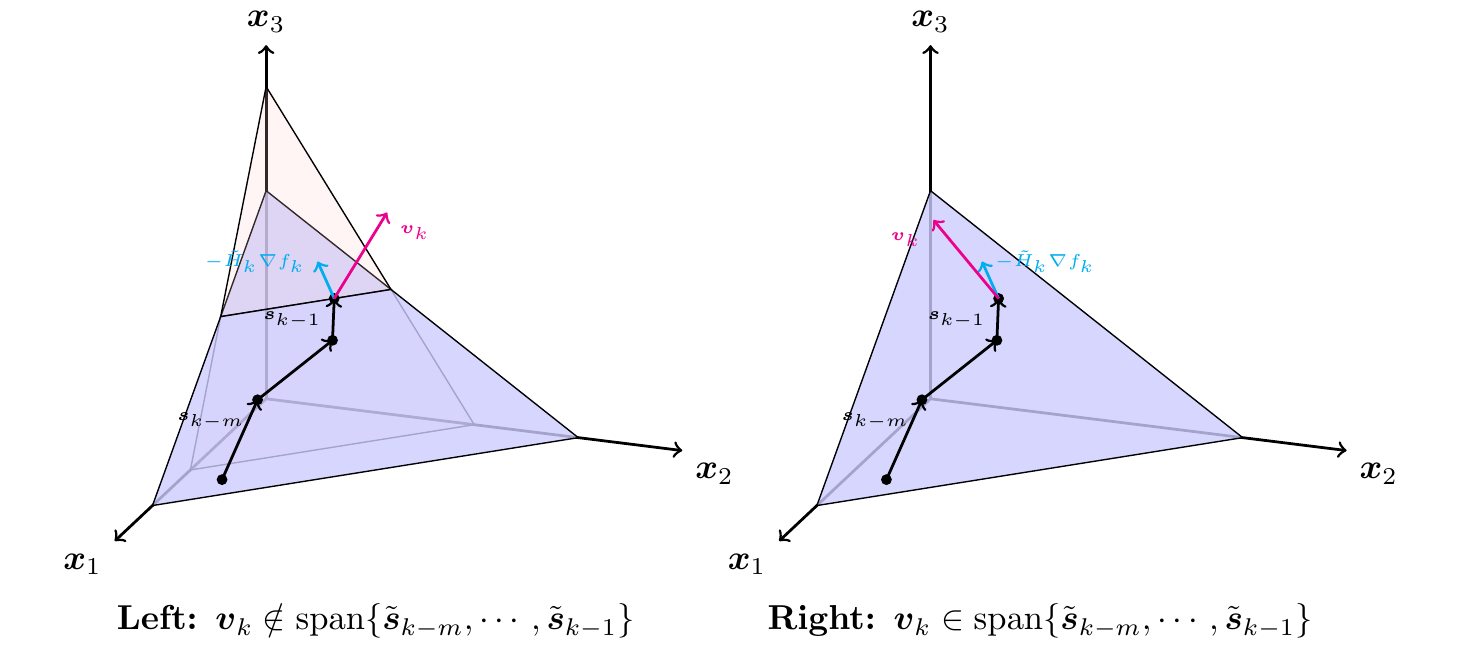}
\caption{\textbf{Left:} $\boldsymbol{v}_k$ helps $\boldsymbol{x}_{k+1}$ to escape from $\boldsymbol{x}_k + \mathrm{span} \{ \tilde{\boldsymbol{s}}_{k-m}, \cdots, \tilde{\boldsymbol{s}}_{k-1} \}$ and $\mathrm{rank} (\tilde{S}_{k+1}) \ge \mathrm{rank} (\tilde{S}_k)$; \textbf{Right:} From Eq. (\ref{eqn: unsolve problem}), $\boldsymbol{p}_k$ is already the best search direction in the local area $B_{\epsilon}(\boldsymbol{x}_k) = \{ \boldsymbol{x} \in \mathbb{R}^n| d(\boldsymbol{x}_k, \boldsymbol{x}) < \epsilon \}$.}
\label{fig: 1}
\end{figure}

\subsection{Estimated the Length of \mathinhead{v}{v}}
\label{sec: Estimated the Length of v}

However the problem (\ref{eqn: unsolve problem}) is too difficult to be solved. If we only consider decaying $\{ \|\nabla f_k\| \}$ with the direction of $\boldsymbol{v}$ is given and length is unknown, the problem becomes easier. Treating $\boldsymbol{v}_k$ as an given unit vector, then modify the search direction as
\begin{equation}
\boldsymbol{p}_k = -\tilde{H}_k \nabla f_k - \alpha_k \boldsymbol{v}_k
\label{eqn: BFGS pk}
\end{equation}
where 
\begin{equation*}
\alpha_k = \mathop{\arg\min}_{\alpha \sim O(\|\nabla f_k\|_2)} \| \nabla f(\boldsymbol{x}_k - \tilde{H}_k \nabla f_k - \alpha \boldsymbol{v}_k) \|_2.
\end{equation*}

Our task now is trun to estimate $\alpha_k$. Defining
\begin{equation*}
J(\alpha) = \|\nabla f(\boldsymbol{x}_k - \tilde{H}_k \nabla f_k - \alpha \boldsymbol{v}_k)\|_2^2,
\end{equation*}
then we obtain
\begin{equation*}
J'(\alpha) = - \boldsymbol{v}^T \nabla^2 f(\boldsymbol{x}_k - \tilde{H}_k \nabla f_k - \alpha \boldsymbol{v}) \nabla f(\boldsymbol{x}_k - \tilde{H}_k \nabla f_k - \alpha \boldsymbol{v}).
\end{equation*}
Let $\epsilon = O(\|\nabla f_k\|_2)$, if $\nabla^3 f(B(\boldsymbol{x}_k; \| \nabla f_k \|_2))$ is bounded, we have
\begin{equation}
J'(\alpha) |_{\alpha \sim \epsilon} = - \boldsymbol{v}_k^T (\nabla^2 f_k + \epsilon) (\nabla f_k - \nabla^2 f_k (\tilde{H}_k \nabla f_k + \alpha \boldsymbol{v}_k) + \epsilon^2).
\label{eqn: gradients estimate}
\end{equation}
Denoting
\begin{equation*}
\alpha_k^* = \frac{\boldsymbol{v}_k^T \nabla^2 f_k (\nabla f_k - \nabla^2 f_k \tilde{\tilde{H}}_k \nabla f_k)}{\boldsymbol{v}_k^T \nabla^2 f_k \nabla^2 f_k \boldsymbol{v}_k}.
\end{equation*}
If $\{\nabla^2 f_k \tilde{H}_k\}_{k \ge m}$ is bounded and $\alpha_k^* \neq 0$, it is simple to show that $\alpha_k^* \sim O(\| \nabla f_k \|_2) = \epsilon$. Applying $\alpha_k$ to Eq. (\ref{eqn: gradients estimate}), we have
\begin{align*}
\boldsymbol{0} = J'(\alpha) |_{\alpha = \alpha_k} 
&= -\boldsymbol{v}_k^T (\nabla^2 f_k + \epsilon) (\nabla f_k - \nabla^2 f_k (\tilde{H}_k \nabla f_k + \alpha_k \boldsymbol{v}_k) + \epsilon^2) \\
&= -\boldsymbol{v}_k^T (\nabla^2 f_k + \epsilon) (\nabla f_k - \nabla^2 f_k \tilde{H}_k \nabla f_k) + \alpha_k \boldsymbol{v}_k^T (\nabla^2 f_k + \epsilon) \nabla^2 f_k \boldsymbol{v}_k + \epsilon^2 \\
&= -\boldsymbol{v}_k^T \nabla^2 f_k (\nabla f_k - \nabla^2 f_k \tilde{H}_k \nabla f_k) + \alpha_k \boldsymbol{v}_k^T \nabla^2 f_k \nabla^2 f_k \boldsymbol{v}_k + \epsilon^2 \\
&= (\boldsymbol{v}_k^T \nabla^2 f_k \nabla^2 f_k \boldsymbol{v}_k) (-\alpha_k^* + \alpha_k) + \epsilon^2.
\end{align*}
Thus
\begin{equation*}
\frac{\alpha_k - \alpha_k^*}{\alpha_k^*} = \frac{-1}{\boldsymbol{v}_k^T \nabla^2 f_k \nabla^2 f_k \boldsymbol{v}_k} \frac{\epsilon^2}{\alpha_k^*} = \epsilon.
\end{equation*}
It follows that 
\begin{equation*}
\lim_{\|\nabla f_k\| \rightarrow 0} \frac{\alpha_k - \alpha_k^*}{\alpha_k^*} = 0.
\end{equation*}
Finally, we estimate $\alpha_k$ by using
\begin{equation}
\alpha_k = \frac{\boldsymbol{v}_k^T \nabla^2 f_k (\nabla f_k - \nabla^2 f_k \tilde{\tilde{H}}_k \nabla f_k)}{\boldsymbol{v}_k^T \nabla^2 f_k \nabla^2 f_k \boldsymbol{v}_k}.
\label{eqn: alpha estimate}
\end{equation}

\subsection{Set the direction of \mathinhead{v}{v}}

In order to make $\{ \| \nabla f_k \| \}$ descend at the beginning of iteration, we should make $\boldsymbol{v}_k$ to satisfy
\begin{equation*}
\alpha_k \boldsymbol{v}_k^T \nabla f_k > 0.
\end{equation*}
Applying Eq. (\ref{eqn: alpha estimate}) to this condition, we get
\begin{equation*}
\boldsymbol{v}_k^T (\nabla^2 f_k \nabla f_k - \nabla^2 f_k \nabla^2 f_k \tilde{H}_k \nabla f_k)  \boldsymbol{v}_k^T \nabla f_k > 0.
\end{equation*}
Thus, if $\nabla^2 f_k \nabla f_k - \nabla^2 f_k \nabla^2 f_k H_k \nabla f_k$ is not parallel to $\nabla f_k$, it easy to show that
\begin{equation*}
\alpha_k (w_1 \boldsymbol{u}_1^{\perp} + w_2 \boldsymbol{u}_2^{\perp})^T \nabla f_k > 0,
\end{equation*}
where $w_1 w_2 > 0$ and
\begin{equation*}
\begin{cases}
\boldsymbol{u}_1 &= \nabla^2 f_k \nabla f_k - \nabla^2 f_k \nabla^2 f_k \tilde{H}_k \nabla f_k; \\
\boldsymbol{u}_2 &=  \nabla f_k; \\
\boldsymbol{u}_1^{\perp} &=  \boldsymbol{u}_1 - \frac{\boldsymbol{u}_1^T \boldsymbol{u}_2}{\|\boldsymbol{u}_1\|_2 \|\boldsymbol{u}_2\|_2} \boldsymbol{u}_2; \\
\boldsymbol{u}_2^{\perp} &=  \boldsymbol{u}_2 - \frac{\boldsymbol{u}_1^T \boldsymbol{u}_2}{\|\boldsymbol{u}_1\|_2 \|\boldsymbol{u}_2\|_2} \boldsymbol{u}_1.
\end{cases}
\end{equation*}
In this paper, we set
\begin{align}
\mathrm{(ver-A)}\quad \boldsymbol{v}_k &= 
\begin{cases}
\left( \frac{\boldsymbol{u}_1^{\perp}}{\|\boldsymbol{u}_1^{\perp}\|_2} + \frac{\boldsymbol{u}_2^{\perp}}{\|\boldsymbol{u}_2^{\perp}\|_2} \right) / \left\| \frac{\boldsymbol{u}_1^{\perp}}{\|\boldsymbol{u}_1^{\perp}\|_2} + \frac{\boldsymbol{u}_2^{\perp}}{\|\boldsymbol{u}_2^{\perp}\|_2} \right\|_2, \quad &\mathrm{if} \ \boldsymbol{u}_1 \neq \boldsymbol{u}_1^{\perp} \ \mathrm{and}\ \boldsymbol{u}_1 \neq \boldsymbol{u}_1^{\perp}; \\
(\boldsymbol{u}_1 + \boldsymbol{u}_2) / \| \boldsymbol{u}_1 + \boldsymbol{u}_2 \|_2, \quad &\mathrm{if} \ \boldsymbol{u}_1^T \boldsymbol{u}_2 = \|\boldsymbol{u}_1\|_2 \|\boldsymbol{u}_2\|_2; \\
\boldsymbol{0}, \quad &\mathrm{if} \ \boldsymbol{u}_1^T \boldsymbol{u}_2 = -\|\boldsymbol{u}_1\|_2 \|\boldsymbol{u}_2\|_2,
\label{eqn: BFGS ver-A vk}
\end{cases} \\
\mathrm{(ver-A)}\quad \alpha_k &= 
\begin{cases}
\frac{\boldsymbol{v}_k^T \nabla^2 f_k (\nabla f_k - \nabla^2 f_k \tilde{H}_k \nabla f_k)}{\boldsymbol{v}_k^T \nabla^2 f_k \nabla^2 f_k \boldsymbol{v}_k}, &\mathrm{if} \ \boldsymbol{u}_1^T \boldsymbol{u}_2 \neq -\|\boldsymbol{u}_1\|_2 \|\boldsymbol{u}_2\|_2; \\
0, \quad &\mathrm{if} \ \boldsymbol{u}_1^T \boldsymbol{u}_2 = -\|\boldsymbol{u}_1\|_2 \|\boldsymbol{u}_2\|_2,
\end{cases}
\label{eqn: BFGS ver-A alphak}
\end{align}
and we don't need to calculate $\nabla^2 f_k$ precisely. Instead of using
\begin{equation}
\begin{cases}
\nabla^2 f_k \nabla^2 f_k \tilde{H}_k \nabla f_k \approx \frac{\nabla f(\boldsymbol{x}_k + \nabla^2 f_k \tilde{H}_k \nabla f_k) - \nabla f_k}{\epsilon}, \quad& \epsilon = \frac{10^{-6}}{\|\nabla^2 f_k \tilde{H}_k \nabla f_k\|_2},\ \mathrm{if}\ \tilde{H}_k \nabla f_k \neq \boldsymbol{0}; \\
\nabla^2 f_k \tilde{H}_k \nabla f_k \approx \frac{\nabla f(\boldsymbol{x}_k + \tilde{H}_k \nabla f_k) - \nabla f_k}{\epsilon}, \quad& \epsilon = \frac{10^{-6}}{\|\tilde{H}_k \nabla f_k\|_2},\ \mathrm{if}\ \tilde{H}_k \nabla f_k \neq \boldsymbol{0} ; \\
\nabla^2 f_k \boldsymbol{v}_k \approx \frac{\nabla f(\boldsymbol{x}_k + \epsilon \boldsymbol{v}_k) - \nabla f_k}{\epsilon}, \quad& \epsilon = \frac{10^{-6}}{\|\boldsymbol{v}_k\|_2}. 
\end{cases}
\label{eqn: BFGS estimate}
\end{equation}
\begin{remark}
If we ignore condition $\alpha_k \boldsymbol{v}_k^T \nabla f_k > 0$ and setting
\begin{align}
\mathrm{(ver-B)}\quad \boldsymbol{v}_k &= \nabla f_k;
\label{eqn: BFGS ver-B vk} \\
\mathrm{(ver-B)}\quad \alpha_k &= \frac{\boldsymbol{v}_k^T \nabla^2 f_k (\nabla f_k - \nabla^2 f_k \tilde{H}_k \nabla f_k)}{\boldsymbol{v}_k^T \nabla^2 f_k \nabla^2 f_k \boldsymbol{v}_k},
\label{eqn: BFGS ver-B alphak}
\end{align}
the sequence $\{ \| f_k \| \}$ generally descending faster, although $\{ f_k \}$ may not be a strictly monotonic-descending sequence in this case. The experiment (Table \ref{tab: CUTE problems}) for CUTE problems illustrating this phenomenon.
\end{remark}

\subsection{Adaptive subspace}
From Section \ref{sec: Estimated the Length of v}, we have $\lim\limits_{\| \nabla f_k | \rightarrow 0} \alpha_k = 0$, it follows that the subspace will tend to gradually stabilize  when the variables converge to the nearby of a minimum point. Thus, the convergence of our algorithm near the minimum point is similar to BFGS.

\subsection{Fast-BFGS Method}

The Complete algorithm of Fast-BFGS is given as

\begin{algorithm}[H]
    \SetAlgoLined
    \KwInput{Small positive number $tol=10^{-5}$.}
    \KwInput{Positive integer $m=8$.}
    \KwInput{The dimension of variables $n$.}
    \KwInput{Initial variables $\boldsymbol{x}_0$.}

    \tcc{Compute $x_1, g_1, \tilde{\boldsymbol{s}}_0, \tilde{\boldsymbol{y}}_0.$}
    Compute search direction $\boldsymbol{p}_0 = -\nabla f_1$; \\
    Set $\boldsymbol{x}_1 = \boldsymbol{x}_0 + \tau_0 \boldsymbol{p}_0$  where $\tau_0$ is computed from a line search procedure to satisfy the strong Wolfe conditions; \\
    Set $\boldsymbol{s}_0 = \boldsymbol{x}_1 - \boldsymbol{x}_0$ and $\boldsymbol{y}_0 = \nabla f_1 - \nabla f_0$; \\
    Define $\tilde{\boldsymbol{s}}_0 = \frac{\boldsymbol{s}_0}{\sqrt{|\boldsymbol{s}_0^T \boldsymbol{y}_0}|}$ and $\tilde{\boldsymbol{y}}_0 = \frac{\boldsymbol{y}_0}{\sqrt{|\boldsymbol{s}_0^T \boldsymbol{y}_0}|}$;

    \tcc{Initial $\tilde{A}_0, \tilde{S}_0, \tilde{S}_1.$}
    Set $\tilde{A}_0 = \begin{bmatrix} 1 \end{bmatrix},
    \tilde{S}_1 = \begin{bmatrix} \tilde{\boldsymbol{s}}_0 \end{bmatrix}, k=1$;

    \tcc{Update $\boldsymbol{x}_{k+1}$.}
    \While{$\|\nabla f_k\| > tol$}
    {
        Compute $\boldsymbol{v}_k$ from Eq. (\ref{eqn: BFGS ver-A vk}) or  (\ref{eqn: BFGS ver-B vk}); \\
        Compute $\alpha_k$ from Eq. (\ref{eqn: BFGS ver-A alphak}) or  (\ref{eqn: BFGS ver-B alphak}); \\
        Compute the search direction $\boldsymbol{p}_k = -\tilde{S}_k \tilde{L}_{k-1} \tilde{S}_k^T \nabla f_k - \alpha_k \boldsymbol{v}_k$; \\
        Set $\boldsymbol{x}_{k+1} = \boldsymbol{x}_k + \tau_k \nabla \boldsymbol{p}_k$ where $\tau_k$ is computed from a line search procedure to satisfy the strong Wolfe conditions; \\
        Set $\boldsymbol{s}_k = \boldsymbol{x}_{k+1} - \boldsymbol{x}_k$ and $\boldsymbol{y}_k = \nabla f_{k+1} - \nabla f_k$; \\
        Compute $\tilde{\boldsymbol{s}}_k = \frac{\boldsymbol{s}_k}{\sqrt{|\boldsymbol{s}_k^T \boldsymbol{y}_k}|}$ and $\tilde{\boldsymbol{y}}_k = \frac{\boldsymbol{y}_k}{\sqrt{|\boldsymbol{s}_k^T \boldsymbol{y}_k|}}$; \\
        Update $\tilde{S}_{k+1}$; \\
        Compute $\tilde{L}_k$ by means of Eq. (\ref{eqn: Fast-BFGS}); \\
        $k = k+1$;
    }
    \KwOutput{$\boldsymbol{x}_{k+1}$.}
    \caption{Parallel in Time BFGS method}
    \label{alg: Fast-BFGS}
\end{algorithm}

\section{Experiments}
\label{sec: Experiments}

We report some numerical results of the Fast-BFGS algorithm. The code is written in Python 3.5 and Tensorflow 1.8.0 on the GPU NVIDIA RTX2080Ti with 1635MHz and 4352 shader cores. All codes(include the serial version on the CPU \& parallel version on the GPU) for Fast-BFGS and related experiments are provided on URL \url{https://github.com/LizhengMathAi/F-BFGS}. All the test functions and their initial values are given by the CUTE collection\cite{bongartz1995cute,lukvsan2010modified,andrei2008unconstrained}. 

\begin{figure}[H]
\centering
\includegraphics[width=1.0\textwidth, angle=0]{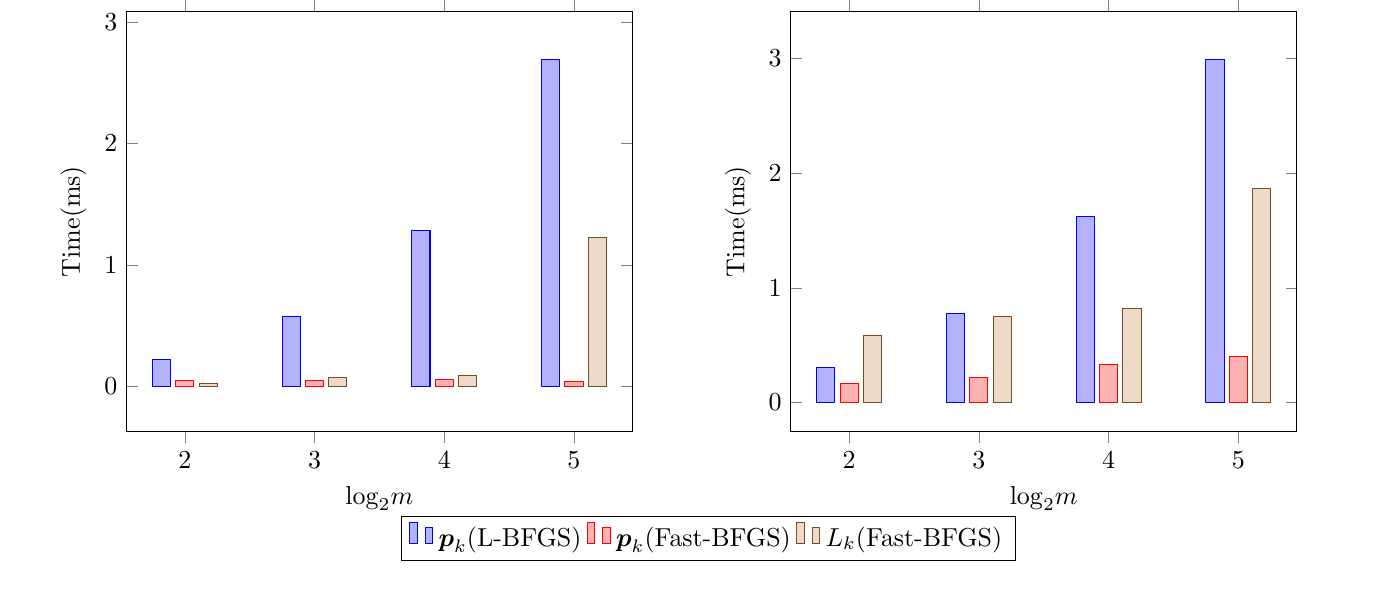}
\caption{\textbf{Left:} $n=2^{10}$; \textbf{Right:} $n=2^{15}$.}
\label{fig: 2}
\end{figure}

When the number of variables is large, the cost of storing $\tilde{H}_{k+1}$ is prohibitive. 
The Fast-BFGS method circumvent this problem. 
Table \ref{tab: storage capacity} compares the BFGS and L-BFGS with Fast-BFGS(ver-A) in terms of storage and computational complexity.
\begin{table}
    \centering
    \caption{ }
    \begin{tabular}{c|c|c|c|c}
        \hline
        \multicolumn{2}{c|}{ } & BFGS & L-BFGS & Fast-BFGS(ver-A) \\
        \hline \hline 
        \multirow{2}{*}{Storage} & items & $H_k$ & $\{\boldsymbol{s}_i, \boldsymbol{y}_i\}_{i=k-m}^{k-1}$ & $\{ \tilde{\boldsymbol{s}}_i \}_{i=k-m}^{k-1}, L_{k-1}$ \\
        \cline{2-5}
        & capacity & $n^2$ & $2mn$ & $mn + m ^2$ \\
        \hline
        Computational & serial & $4n^3+o(n^3)$ & $11mn + o(mn)$ & $n(2m^2+8m+23) + o(n)$ \\
        \cline{2-5}
        complexity & parallel & $3n+o(n)$ & $2mn + o(mn)$ & $8n + o(n)$ \\
        \hline
    \end{tabular}
    \label{tab: storage capacity}
\end{table}
Fig. \ref{fig: 2} shows the GPU performance profile of Fast-BFGS versus L-BFGS. We see that Fast-BFGS was better in large-scale parallel computing. 

In the next set of experiments, Table \ref{tab: CUTE problems} results illustrating the behavior of Fast-BFGS and other methods for those large-scale unconstrained optimization problems taken from the CUTE collection. It gives the number of function and gradient evaluations (nfg) and the termination criterion $\| \nabla f_k \|_2 < 10^{-5}$ is used. 
\begin{table}[H]
    \centering
    \caption{Performance of Fast-BFGS method.}
    \label{tab: CUTE problems}
    \begin{tabular}{|lr|c|c|c|c|c|}
    \hline
             &      & GD    & BFGS  & L-BFGS & ver-A & ver-B \\ \cline{3-7}
    Problem  & n    & \multicolumn{2}{c|}{nfg} & \multicolumn{3}{c|}{nfg(m=8)} \\ \hline\hline
    ARWHEAD  & 1024 & >1000 & 39    & 26     & 21    & 16    \\ \hline
    BDQRTIC  & 1024 & >1000 & --    & --     & 491   & 317   \\ \hline
    BDEXP    & 1024 & >1000 & 19    & 19     & 9     & 9     \\ \hline
    COSINE   & 1024 & >1000 & --    & --     & 44    & 16    \\ \hline
    DIXMAANE & 1500 & >1000 & 195   & 244    & 586   & 326   \\ \hline
    DIXMAANF & 1500 & >1000 & 336   & 216    & 423   & 265   \\ \hline
    DIXMAANG & 1500 & >1000 & 954   & 384    & 460   & 211   \\ \hline
    DQRTIC   & 1000 & --    & --    & --     & 35    & 31    \\ \hline
    EDENSCH  & 1000 & 59    & 86    & 52     & 42    & 23    \\ \hline
    ENGVAL1  & 1000 & 66    & 154   & 119    & 39    & 24    \\ \hline
    EG2      & 1000 & 7     & 6     & 6      & 8     & 8     \\ \hline
    EXTROSNB & 1000 & 63    & 309   & 333    & 76    & 41    \\ \hline
    FLETCHER & 100  & >1000 & --    & --     & >1000 & 734   \\ \hline
    FREUROTH & 1000 & --    & --    & --     & 51    & 45    \\ \hline
    GENROSE  & 1000 & >1000 & >1000 & 39     & 48    & --    \\ \hline
    HIMMELBG & 1000 & >1000 & 3     & 3      & 3     & 3     \\ \hline
    HIMMELH  & 1000 & 20    & 9     & 9      & 19    & 16    \\ \hline
    LIARWHD  & 1000 & >1000 & --    & 28     & 40    & 30    \\ \hline
    NONDIA   & 1000 & >1000 & --    & 55     & 97    & 76    \\ \hline
    NONDQUAR & 1000 & >1000 & 270   & 320    & 344   & 230   \\ \hline
    NONSCOMP & 1000 & 86    & 286   & 238    & 101   & 45    \\ \hline
    POWELLSG & 1000 & >1000 & 459   & 49     & 69    & 63    \\ \hline
    SCHMVETT & 1000 & 181   & 26    & 24     & 45    & 25    \\ \hline
    SINQUAD  & 1000 & >1000 & 140   & 143    & --    & --    \\ \hline
    SROSENBR & 1000 & >1000 & --    & 39     & 48    & --    \\ \hline
    TOINTGSS & 1000 & 6     & 9     & 9      & 8     & 7     \\ \hline
    TQUARTIC & 1000 & >1000 & 16    & 17     & 28    & 24    \\ \hline
    WOODS    & 1000 & >1000 & --    & 92     & --    & 48    \\ \hline
    \end{tabular}
\end{table}

Table \ref{tab: CUTE problems} shows that our algorithm is more effective than BFGS method or L-BFGS method most of the time. Our conjecture for this phenomenon is, 
the inverse Hessian $\left( \nabla^2 f_k \right)^{-1}$ may sometimes be ill-conditioned so that BFGS and L-BFGS methods cannot approximate the $\left( \nabla^2 f_k \right)^{-1} \nabla f_k$. But our algorithm is actually estimating the dynamic inverse Hessian matrix $\left( \nabla_{\boldsymbol{\xi}}^2 f_k \right)^{-1}$ on a low-dimensional space $\mathbb{R}^m$, which is less likely to be ill-conditioned.

In the end, Table \ref{tab: weakly dependent} presents results illustrating the behavior of Fast-BFGS for various levels of memory $m$, it shows that our method tends to be still robust when $m$ is small at most of the time.
\begin{table}[H]
\centering
\caption{The rate of convergence is weakly dependent on $m$.}
\label{tab: weakly dependent}
\begin{tabular}{|lr|cc|cc|cc|}
\hline
         &      & ver-A & ver-B & ver-A & ver-B & ver-A & ver-B \\ \cline{3-8}
Problem  & n    & \multicolumn{2}{c|}{nfg(m=2)} & \multicolumn{2}{c|}{nfg(m=4)} & \multicolumn{2}{c|}{nfg(m=8)} \\ \hline\hline
ARWHEAD  & 1024 & 21    & 16    & 21    & 16    & 21    & 16    \\ \hline
BDQRTIC  & 1024 & >1000 & >1000 & >1000 & 427   & 491   & 317   \\ \hline
BDEXP    & 1024 & 9     & 9     & 9     & 9     & 9     & 9     \\ \hline
COSINE   & 1024 & 98    & 17    & 63    & 16    & 44    & 16    \\ \hline
DIXMAANE & 1500 & 586   & 800   & 619   & 327   & 586   & 326   \\ \hline
DIXMAANF & 1500 & 513   & 559   & 481   & 252   & 423   & 265   \\ \hline
DIXMAANG & 1500 & 520   & 233   & 457   & 268   & 460   & 211   \\ \hline
DQRTIC   & 1000 & 36    & 32    & 35    & 31    & 35    & 31    \\ \hline
EDENSCH  & 1000 & 49    & 31    & 46    & 28    & 42    & 23    \\ \hline
ENGVAL1  & 1000 & 55    & 30    & 44    & 26    & 39    & 24    \\ \hline
EG2      & 1000 & >1000 & >1000 & 8     & 7     & 8     & 8     \\ \hline
EXTROSNB & 1000 & 77    & 44    & 77    & 42    & 76    & 41    \\ \hline
FLETCHER & 100  & >1000 & >1000 & --    & 791   & >1000 & 734   \\ \hline
FREUROTH & 1000 & 248   & 70    & 65    & --    & 51    & 45    \\ \hline
GENROSE  & 1000 & 48    & 55    & 50    & --    & 48    & --    \\ \hline
HIMMELBG & 1000 & 3     & 3     & 3     & 3     & 3     & 3     \\ \hline
HIMMELH  & 1000 & 19    & 16    & 19    & 16    & 19    & 16    \\ \hline
LIARWHD  & 1000 & 40    & 30    & 39    & 30    & 40    & 30    \\ \hline
NONDIA   & 1000 & 90    & 74    & 93    & 74    & 97    & 76    \\ \hline
NONDQUAR & 1000 & 795   & 953   & 571   & 278   & 344   & 230   \\ \hline
NONSCOMP & 1000 & 108   & 58    & 107   & 49    & 101   & 45    \\ \hline
POWELLSG & 1000 & >1000 & 497   & 66    & 63    & 69    & 63    \\ \hline
SCHMVETT & 1000 & 88    & 52    & 79    & 24    & 45    & 25    \\ \hline
SINQUAD  & 1000 & >1000 & 316   & --    & >1000 & --    & --    \\ \hline
SROSENBR & 1000 & 48    & 70    & 48    & 85    & 48    & --    \\ \hline
TOINTGSS & 1000 & 10    & 7     & 8     & 7     & 8     & 7     \\ \hline
TQUARTIC & 1000 & 28    & 23    & 28    & 23    & 28    & 24    \\ \hline
WOODS    & 1000 & 638   & 254   & 48    & 48    & --    & 48    \\ \hline
\end{tabular}
\end{table}

\bibliographystyle{unsrt}  
\bibliography{references} 

\end{document}